# PERFECT SAMPLING USING BOUNDING CHAINS

By Mark Huber[1]

*Duke University*

Bounding chains are a technique that offers three benefits to Markov chain practitioners: a theoretical bound on the mixing time of the chain under restricted conditions, experimental bounds on the mixing time of the chain that are provably accurate and construction of perfect sampling algorithms when used in conjunction with protocols such as coupling from the past. Perfect sampling algorithms generate variates exactly from the target distribution without the need to know the mixing time of a Markov chain at all. We present here the basic theory and use of bounding chains for several chains from the literature, analyzing the running time when possible. We present bounding chains for the transposition chain on permutations, the hard core gas model, proper colorings of a graph, the antiferromagnetic Potts model and sink free orientations of a graph.

**1. Introduction.** The breadth of applications using Monte Carlo Markov chain (MCMC) techniques today is astounding. From theoretical physics to approximation of $\sharp P$-complete problems to Bayesian inference, all of these techniques rest on the ability to construct a Markov chain whose stationary distribution matches a particular target distribution whose normalizing constant is unknown. Many methods (such as Gibbs sampling and Metropolis–Hastings [15, 22]) exist for constructing such chains, but finding the mixing time of these chains is, in general, extremely difficult. Although some theoretical progress on this problem has been made [5], for most chains of practical interest, heuristics, such as autocorrelation tests [11], are the only methods available. Such methods only give negative results—they can indicate when a Markov chain has not mixed, but do not guarantee that the

Received August 2002; revised December 2002.
[1]Supported by an ONR Fellowship, NSF Grants CCR-93-07391, CCR-97-00029 and DMS-95-05155, ONR Grant N00014-96-1-00500 and NSF postdoc 99-71064.

random variates obtained come from the correct distribution or anything remotely close to it. Estimates based upon such samples will necessarily be suspect.

Bounding chains, when applicable, cut through this problem and give true algorithms for obtaining random variates. They are useful for MCMC in three ways. First, they are formed from a natural extension of couplings we shall refer to as complete couplings. When it is possible to analyze the bounding chain theoretically, this relationship to couplings allows derivation of a theoretical bound on the mixing time of the Markov chain based on well-known results.

Second, bounding chains are themselves Markov chains that can be simulated. Again because of the relationship to couplings, simulations of the bounding chain can give insight into the probability of coupling in the original chain and so experimental evidence can be obtained about the mixing time of the original chain. To be clear, this evidence is very different from that obtained from autocorrelation tests. Autocorrelation tests can only give evidence that a chain has not mixed, and the more complex the chain, the less value such tests have. Use of bounding chains can give direct evidence that a chain is mixing rapidly and the quality of this evidence can be analyzed precisely.

The third use is the most powerful: bounding chains allow a user to utilize perfect sampling algorithms such as the coupling from the past (CFTP) method of Propp and Wilson [23] or the method of Fill, Machida, Murdoch and Rosenthal [10]. Perfect sampling techniques (also known as perfect simulation methods) draw samples exactly from distributions where the normalizing constant of the distribution is unknown. Although both CFTP and [10] utilize Markov chains in their operation, they are true algorithms: there is no need to have an understanding of the mixing time of the underlying Markov chain in order to run these algorithms. At termination the variates generated come exactly from the target distribution.

In this work we will present the basic framework for bounding chains, along with the primary method for their construction. This will then be applied to an introductory example, followed by more advanced applications drawn from statistical mechanics, graph theory and the approximation of $\sharp P$-complete problems.

**2. Bounding chains.** Bounding chains were introduced independently by the author [16] and Häggström and Nelander [14]. Here we present two basic ways of looking at bounding chains that arise naturally from consideration of couplings.

Let $\mathcal{M}$ be a Markov chain with state space $\Omega$, kernel $\mathbf{K}$ and stationary distribution $\pi$. Our goal is to generate random variates from $\pi$. The classical method of starting the Markov chain at a state $x_0$ and running for a large number of steps requires some means for measuring how close the distribution is to stationarity, such as the *total variation distance*: $\|\mathbf{K}^t(x_0, \cdot) - \pi(\cdot)\|_{\text{TV}} = \sup_A |\mathbf{K}^t(x_0, A) - \pi(A)|$.

Let $\mathcal{M}_1$ and $\mathcal{M}_2$ be two chains with kernels $\mathbf{K}_1$ and $\mathbf{K}_2$ and state spaces $\Omega_1$ and $\Omega_2$. A coupling between $\mathcal{M}_1$ and $\mathcal{M}_2$ is a bivariate process $(X_t, Y_t)$ on $\Omega_1 \times \Omega_2$, where the marginal distribution of $X_t$ is Markov with kernel $\mathbf{K}_1$ and $Y_t$ is also marginally Markov with kernel $\mathbf{K}_2$.

Suppose we couple a chain with a copy of itself to give the bivariate process $(X_t, Y_t)$. We shall refer to this particular coupling as a simple, or pairwise, coupling. Let $X_0 = x_0$, $Y_0 \sim \pi$, and suppose that $X_t = Y_t \Rightarrow X_{t'} = Y_{t'}$ for all $t' \geq t$. Under this (very strong) condition we may employ the *Coupling lemma* [1, 7] to state:

$$\|\mathbf{K}^t(x_0, \cdot) - \pi(\cdot)\|_{\text{TV}} \leq P(X_t \neq Y_t).$$

Perfect sampling protocols such as coupling from the past require a more global form of coupling. Suppose that we have a family of measurable functions from $\Omega$ to $\Omega$ and a probability distribution on the family such that for a draw $\phi$, $P(\phi(x) \in A) = \mathbf{K}(x, A)$ for all $x \in \Omega$ and measurable $A$. If $\phi_0, \phi_1, \ldots$ are a sequence of independent draws from this family of distributions, then for any $x_0 \in \Omega$, $x_0, \phi_0(x_0), \phi_1(\phi_0(x_0)), \ldots$ will be a simulation of the Markov chain starting from state $x_0$. Of course, if we start from $x_0$ and $y_0$ using the same draws $\phi_0, \phi_1, \ldots$, then we have a simple coupling satisfying the conditions of the Coupling lemma. In this paper we shall refer to such a family of functions and its distribution as a *complete coupling*. Let $\Phi_0$ be the identity function and $\Phi_t = \phi_{t-1} \circ \Phi_{t-1}$ for $t \geq 1$. Then if $\Phi_t(\Omega) = \{c\}$, that is, if $\Phi_t$ is a constant function, then every single state in the state space has coupled. This immediately gives us:

$$\|\mathbf{K}^t(x_0, \cdot) - \pi(\cdot)\|_{\text{TV}} \leq P(\Phi_t \text{ is not constant}).$$

CFTP and related methods generate variates drawn exactly from the target distribution $\pi$. Three elements are required: a Markov chain with $\pi$ as its stationary distribution, a complete coupling for the chain and a means of determining when $\Phi_t$ is constant. CFTP computes $\Phi_t$ for a fixed value of $t$. If $\Phi_t$ is a constant state $x$ after these $t$ steps, then $x$ is output as the random variate. If $\Phi_t$ is not constant, CFTP is called recursively to generate a random variate $y$. Then $x = \Phi_t(y)$ is output as the random variate. As long as there is a reasonably large probability $\Phi_t$ is a constant, the number of recursive calls within CFTP will be small. The output of this algorithm is distributed exactly according to $\pi$, the stationary distribution of the chain [23].

One effective means for determining complete coupling in order to use CFTP is to take advantage of a monotonic structure in the chain. We will not go into the details here, other than to note that this idea was later extended to a wider variety of chains by introducing antimonotonicity [13]. Still, there remain chains of interest that are neither monotone nor antimonotone.

Bounding chains are a method for determining when $\Phi_t$ is constant for a wide variety of state spaces. The two forms of bounding chains we present are identical when dealing with discrete state spaces, but act differently when extended to continuous state spaces.

Suppose that the chain $\mathcal{M}$ with state space $\Omega \subseteq C^V$ for some set of dimensions $V$ and some set of possible values $C$ (we shall refer to $C$ as the set of colors). Usually $|\Omega|$ is growing exponentially with $|V|$, which is why MCMC methods are commonly used in this situation. Consider a new Markov chain $\mathcal{M}'$ with state space $(2^C)^V$, where $2^C$ is the set of subsets of $C$. Let $X_t$ be a stochastic process evolving according to $\mathcal{M}$ and $Y_t$ be a stochastic process evolving according to $\mathcal{M}'$.

DEFINITION 1 (Form 1). We will say that $\mathcal{M}'$ is a *bounding chain* for $\mathcal{M}$ if there exists a coupling between $\mathcal{M}'$ and $\mathcal{M}$ such that

$$X_t(v) \in Y_t(v) \quad \forall v \implies X_{t+1}(v) \in Y_{t+1}(v) \quad \forall v.$$

Each dimension $v \in V$ is given a single color from $C$ in $X_t$ and a subset of colors from $C$ in $Y_t$. If we have $\mathcal{M}'$ bounding $\mathcal{M}$, we have a coupling for $(X_t, Y_t)$ so that by starting with $X_0(v) \in Y_0(v)$ for all $v$, we guarantee that $X_t(v) \in Y_t(v)$ for all $v \in V$ and all times $t \geq 0$. If $X_t$ is evolving using $\phi_0, \phi_1, \ldots$, and $Y_0$ bounds every state in $\Omega$, then when $Y_t$ bounds just one state $x$, that means that $\Phi_t$ is a constant function. The number of states bounded by $Y_t$ is at most $\prod_v |Y_t(v)|$, so to check that $\Phi_t$ is constant, it suffices that $|Y_t(v)| = 1$ for all $v$.

Given a complete coupling for $\mathcal{M}$, generated by $\phi_0, \phi_1, \ldots$, we can construct a bounding chain as follows. Say that $x \in C^V$ is bounded by $y \in (2^C)^V$ (denote $x \in_V y$) if $x(v) \in y(v)$ for all $v$. Then given state $(X_t, Y_t)$, let $X_{t+1} = \phi_t(X_t)$ and let

$$Y_{t+1}(v) = \bigcup_{x \in_V Y_t} (\phi_t(x))(v).$$

That is, consider each and every possible state $x$ that is bounded by $Y_t$. Take the next step from each of these states $x$ and examine the set of colors for $v$ that results. By taking $Y_{t+1}(v)$ to be this set of colors and doing this for every node $v \in V$, we guarantee that $Y_{t+1}$ bounds the state $X_{t+1}$. Starting the bounding chain in a state $Y_0$ that bounds $X_0$ is easy when $V$ and $C$ are discrete: just make $Y_0(v) = C$ for all $v \in V$. Examples of problems where Form 1 is the most useful formulation include proper colorings of a graph and the Potts model.

Sometimes it is more convenient to write the state space as $\Omega = C^{(2^V)}$, so if we think of a configuration as a function mapping $V$ to $C$, instead keep track of the inverse of the function for each $c \in C$. That is, for each color, we record the subset of dimensions that have that color. In this case, the bounding chain $\mathcal{M}''$ will have state space $\{(b,d) : b, d \in \Omega, b(c) \cap d(c) = \varnothing \ \forall c\}$, that is, ordered pairs of states from $\Omega$ that do not intersect. Let $A_t$ be a stochastic process evolving according to $\mathcal{M}$ and $(B_t, D_t)$ a stochastic process evolving according to $\mathcal{M}''$.

DEFINITION 2 (Form 2). We will say that $\mathcal{M}''$ is a *bounding chain* for $\mathcal{M}$ if there exists a coupling between $\mathcal{M}''$ and $\mathcal{M}$ such that

$$B_t(c) \subseteq A_t(c) \subseteq B_t(c) \cup D_t(c) \qquad \forall c$$
$$\implies B_{t+1}(c) \subseteq A_{t+1}(c) \subseteq B_{t+1}(c) \cup D_{t+1}(c) \qquad \forall c.$$

Intuitively, for $A_t$ bounded by $(B_t, D_t)$, $B_t(c)$ is the set of dimensions that are "known" to have color $c$ in $A_t$, while $D_t$ are those dimensions that "might" have color $c$ in $A_t$. As with the other form, complete coupling has occurred when $D_t$ is empty, which is very easy to check.

We can construct a bounding chain in Form 2 using a similar method to the first form. Say that $a$ is bounded by $(B_t, C_t)$ if for every $c \in C$ we have $B_t(c) \subseteq a \subseteq B_t(c) \cup D_t(c)$ and we write $a \in_C (B_t, D_t)$. Given $\phi_0, \phi_1, \ldots$ for a complete coupling, let $A_{t+1} = \phi_t(A_t)$, and for each $c \in C$, set

$$B_{t+1}(c) = \bigcap_{a \in_C (B_t, D_t)} (\phi_t(a))(c),$$

$$D_{t+1}(c) = \left(\bigcup_{a \in_C (B_t, D_t)} (\phi_t(a))(c)\right) \setminus B_{t+1}.$$

Start the bounding chain in state $(B_0, D_0)$, where $B_0(c) = \varnothing$ and $D_0(c) = V$ for all $c \in C$. This will trivially bound $A_0$.

Examples of problems where Form 2 is more useful include such chains as Swendsen–Wang [17] and independent sets of a graph [19].

**3. Transposition chain for permutations.** Our first application of the bounding chain technique is a toy example intended to illustrate several aspects of the method. While no practitioner would use this chain for generation of random permutations, the techniques we use here will be applicable to more sophisticated examples later. Consider the problem of generating random permutations of $n$ objects uniformly at random. Using our coloring scheme, $\Omega$ is the subset of $\{1, \ldots, n\}^{\{1, \ldots, n\}}$ such that $x(v) \neq x(w) \ \forall v, w$. If $x(i) = j$, we will say that *item $j$ is in position $i$*. The Markov chain takes steps by randomly moving between states that differ by a single transposition. Specifically we set

$$\mathbf{K}(x,y) = \begin{cases} 2/n^2, & x,y \text{ differ in two positions,} \\ 1/n, & x = y, \\ 0, & \text{otherwise.} \end{cases}$$

There are of course several ways to create complete couplings for this kernel. For instance, $i$ and $j$ could be chosen independently and uniformly at random from $\{1,\ldots,n\}$, and then set $x_{t+1}(i) = x_t(j)$ and $x_{t+1}(j) = x_t(i)$. This is equivalent to switching the items at position $i$ and $j$.

Alternately, items $i$ and $j$ could be swapped, so $x_{t+1}(x_t^{-1}(i)) = j$ and $x_{t+1}(x_t^{-1}(j)) = i$. Although these satisfy the definition of complete couplings, $\Phi_t$ will never become constant for either of these approaches and the bounding chain derived from them will never detect complete coupling.

A more useful complete coupling swaps item $i$ with position $j$. To turn this complete coupling into a bounding chain, we just follow the procedure outlined in the previous section. We will utilize Form 1 (1) in this example. We begin with all the positions assigned state $\{1,\ldots,n\}$ in the bounding chain. As we run the bounding chain, some positions will move to singleton sets of colors, but some will have more than one color. Call a position *known* if it is assigned a set of size 1 in the bounding chain state and call it *unknown* otherwise. Similarly, call an item known if some position is known and assigned that item, and call it unknown otherwise. To minimize the needed updating, each unknown position in the bounding chain will always be assigned $\{1,\ldots,n\}$ even though we could shrink the size of this set somewhat by considering the values of the known positions.

Suppose our bounding chain is in state $Y_t$ and we choose to swap item $i$ with position $j$. Then for any $x$ bounded by $Y_t$, $\phi_t(x)(j) = i$ because we move item $i$ to position $j$. Taking the union of $\{i\}$ over all $x$ bounded by $Y_t$ just gives $Y_{t+1}(j) = \{i\}$.

To find $Y_{t+1}(i)$, we must consider two cases. In the first case, $i$ is known to have position $i_{\text{pos}}$ in $Y_t$ [so $Y_t(i_{\text{pos}}) = \{i\}$]. In this case, $\phi_t(x)(i_{\text{pos}}) = x(j)$. Now the union over all $x$ bounded by $Y_t$ of $x(j)$ is a subset of $Y_t(j)$, so setting $Y_{t+1}(i_{\text{pos}}) = Y_t(j)$ insures that $Y_{t+1}$ bounds $\phi_t(x)$.

In the second case, $i$ has unknown position in $Y_t$, so that the location of $i$ in an $x$ bounded by $Y_t$ is somewhere in the unknown positions in $Y_t$. Hence the item at position $j$ moves to one of these unknown positions at the next time step. But any unknown position $k$ has $Y_t(k) = \{1,\ldots,n\}$ and so setting $Y_{t+1}(k) = Y_{t+1}(k)$ for all $k \neq j$ makes sure that for $x$ bounded by $Y_t$, $\phi_t(x)$ is bounded by $Y_{t+1}$.

Now we analyze the expected time needed for this bounding chain to detect complete coupling, that is, for all of the positions to become known.

LEMMA 1. *Let $\tau$ be the first time the bounding chain detects complete coupling. Then*

(1) $$E[\tau] \leq \frac{\pi^2}{6} n^2 \quad \text{and} \quad P(\tau > CE[\tau]) \leq e^{-0.18C}$$

*for all $C \geq 2$.*

PROOF. Let $W_t$ denote the number of unknown positions in the bounding chain at time $t$. We are interested in computing the expected time for $W_t$ to hit 0. Suppose that item $i$ and position $j$ are selected for the switch. Item $i$ and position $j$ can each be either known or unknown. If at least one of $i$ or $j$ is unknown, then $W_t$ is unchanged. If however, both $i$ and $j$ are unknown, then after one step item $i$ (or equivalently position $j$) becomes known, $W_t$ decreases by $i$. Since $i$ and $j$ were chosen uniformly and independently, the probability that both are unknown is just $(W_t/n)^2$. Hence the expected time until $W_t$ decreases by 1 is just $(n/W_t)^2$ and the expected time until $W_t = 0$ is

(2) $$\sum_{i=0}^{n} \left(\frac{n}{W_t}\right)^2 \leq n^2 \sum_{i=0}^{n} \frac{1}{W_t^2} \leq n^2 \frac{\pi^2}{6}.$$

For the tail bound, note that the probability that $T \geq 2E[\tau]$ is at most $1/2$ by Markov's inequality. So for each $2E[\tau]$ time steps we have at least a $1/2$ chance of detecting complete coupling, so after $CE[\tau]$ steps, the probability that we have failed is at most $(1/2)^{(C/2)}$ for integer $C$. Accounting for the fact that $C$ might be real and not an integer reduces our bound to $(1/2)^{(C/4)} \leq \exp(-0.18C)$. □

In this case, the expected time for the bounding chain to detect complete coupling comes within a constant factor of the expected time needed for any complete coupling from the class we are considering to couple even two processes. Suppose that $\{X_t\}$ and $\{Y_t\}$ are both evolving according to $\phi_0, \phi_1, \ldots$ and let $\tau$ be the amount of time needed for coupling to occur. Note that $X_t$ and $Y_t$ cannot differ in exactly one position and each step changes at most 4 positions. At some point before coupling occurs they differ in 2, 3 or 4 positions. To finish the coupling, either $X_t$ or $Y_t$ must change states and one of $X_t$ and $Y_t$ must choose two positions where they differ. This occurs with probability at most $24/n^2$ and so the expected time to couple will be at least $n^2/24$, differing from the bounding chain upper bound by a constant factor.

While $\Theta(n^2)$ time is needed for this type of coupling, the chain itself is known to mix in $\Theta(n \ln n)$ time [6]. These special couplings do not always tightly bound the mixing time of a Markov chain, a weakness inherited by bounding chains.

**4. The hard core gas model.** An *independent set* of a graph is a subset of nodes such that no two nodes in the subset are connected by an edge of the graph. In this section we wish to generate random variates from the distribution

$$\pi(A) = \frac{\lambda^{|A|}}{Z_\lambda}, \tag{3}$$

where $A$ is any independent set of the graph, $\lambda$ is a parameter of the model known as the *fugacity* and $Z_\lambda$ is the (unknown) normalizing constant often called the *partition function* of the model.

In statistical physics, this distribution is known as the hard core gas model [25] and is a simple model of gases where each node of the graph might or might not be occupied by a gas molecule of nonnegligible size so that any two nodes connected by an edge cannot both be occupied simultaneously. The fugacity controls the density of the gas. Despite its simplicity, the model exhibits a phase transition when the graph satisfies symmetry properties. For instance, the lattice $Z^d$ satisfies these properties. The same distribution also arises from certain stochastic models of communication networks [21].

Computation of the partition function for general graphs is a $\sharp P$-complete problem [8]. Because this is an example of a self-reducible problem, the ability to efficiently generate variates from this distribution immediately leads to a polynomial time approximation algorithm for finding $Z_\lambda$ [20, 24].

In [9], Dyer and Greenhill proposed a Markov chain for this distribution that was neither the straightforward Gibbs sampler nor a Metropolis–Hastings algorithm, but instead a combination of the two that we shall refer to as the Dyer–Greenhill chain. They used path coupling to show that when $\lambda < 2/(\Delta - 2)$, where $\Delta$ is the maximum degree of the graph, the chain mixes in $\Theta(n \ln n)$ time. Their chain is neither monotonic nor antimonotonic. We present a complete coupling of their chain, develop the corresponding bounding chain and show the following:

THEOREM 1. *When $\lambda < 2/(\Delta - 2)$, the bounding chain for the Dyer–Greenhill chain detects complete coupling after $[\log_\beta n] + \theta$ steps with probability at least $1 - \beta^\theta$, where*

$$\beta = \frac{\Delta \lambda}{2(\lambda + 1)}. \tag{4}$$

The state of the chain at time $t$ can be described by a single subset of nodes $A_t$ that is the independent set. Using Form 2 gives us a bounding chain whose state is a single ordered pair $(B_t, D_t)$. A regular Gibbs sampler chooses a node at random, and then adds or deletes that node with probabilities that depend on $\lambda$ and whether neighbors of the node are in $A_t$. The Dyer–Greenhill chain adds a new wrinkle: if exactly one neighbor of the chosen vertex is in the independent set, the new node can "swap," where it is added to the independent set and the neighbor is deleted. This swap move occurs with probability $p_{\text{swap}}$, a new parameter of the chain. Let $s \in_U S$ mean that we choose an element $s$ uniformly at random from the set $S$. Given a graph with node set $\{1, \ldots, n\}$, the following is a complete coupling for their chain that takes a single step from state $A$:

*The Dyer–Greenhill chain for the hard core gas model*
1.   **Choose** $v \in_U \{1, \ldots, n\}$
2.   **Choose** $U \in_U [0,1]$
3.   **Case I:** $U > \lambda/(\lambda+1)$
4.     Let $A \leftarrow A \setminus \{v\}$
5.   **Case II:** $U < \lambda/(\lambda+1)$, and no neighbors of $v$ are in $A$
6.     Let $A \leftarrow A \cup \{v\}$
7.   **Case III:** $U < p_{\text{swap}}\lambda/(\lambda+1)$, exactly one neighbor of $v$ (call it $w$) is in $A$
8.     Let $A \leftarrow A \setminus \{w\} \cup \{v\}$

Lines 1 and 2 are the random choice of $\phi$ for this time step and the remaining lines merely compute $\phi(A)$. For the bounding chain we choose $v$ and $U$ in the same way and then see how $(B, D)$ changes based on those decisions. Using Form 2 (2), we wish to insure that $B \subseteq A \subseteq B \cup D$.

The bounding chain considers six cases. Case I of the original chain corresponds to Case I below, where the state of the neighbors of $v$ do not matter at all because we have rolled to remove $v$ from the independent set. Case II corresponds to cases in the chain where no neighbors of $v$ are in $B$ and in Case III exactly one neighbor of $v$ will be in $B$.

*Bounding chain for the Dyer–Greenhill chain*
1. **Choose** vertex $v \in_U \{1, \ldots, n\}$, let $N_v$ be the neighbors of $v$
2. **Choose** $U \in_U [0,1]$
3. **Case I:** $U > \lambda/(\lambda+1)$
4.   Let $B \leftarrow B \setminus \{v\}$, $D \leftarrow D \setminus \{v\}$,
5. **Case IIa:** $U < \lambda/(\lambda+1)$, $|N_v \cap B| = |N_v \cap D| = 0$
6.   Let $B \leftarrow B \cup \{v\}$, $D \leftarrow D \setminus \{v\}$
7. **Case IIb:** $p_{\text{swap}}\lambda/(\lambda+1) \leq U < \lambda/(\lambda+1)$, $|N_v \cap B| = 0$, $|N_v \cap D| = 1$
8.   Let $B \leftarrow B \setminus \{v\}$, $D \leftarrow D \cup \{v\}$
9. **Case IIc:** $U < p_{\text{swap}}\lambda/(\lambda+1)$, $|N_v \cap B| = 0$, $|N_v \cap D| = 1$
10.   Let $B \leftarrow B \cup \{v\}$, $D \leftarrow D \setminus \{v\} \setminus (N_v \cap D)$
11. **Case IId:** $U < \lambda/(\lambda+1)$, $|N_v \cap B| = 0$, $|N_v \cap D| \geq 2$
12.   Let $B \leftarrow B \setminus \{v\}$, $D \leftarrow D \cup \{v\}$
13. **Case IIIa:** $U < p_{\text{swap}}\lambda/(\lambda+1)$, $|N_v \cap B| = 1$, $|N_v \cap D| = 0$
14.   Let $B \leftarrow B \cup \{v\} \setminus (N_v \cap B)$, $D \leftarrow D \setminus \{v\}$
15. **Case IIIb:** $U < p_{\text{swap}}\lambda/(\lambda+1)$, $|N_v \cap B| = 1$, $|N_v \cap D| \geq 1$
16.   Let $B \leftarrow B \setminus \{v\} \setminus (N_v \cap B)$, $D \leftarrow D \cup \{v\} \cup (N_v \cap B)$

Note that Case II has four subcases. If no neighbors of $v$ are in $D$ as well, then the node is always added for any $A$ bounded by $(B, D)$. If two or more neighbors of $v$ are in $D$, then $v$ might be added or it might not, so we must add it to $D$.

If exactly one neighbor of $v$ is in $D$, then if we do not swap based on $U$, we could end up with $v$ in the independent set or not for a state bounded by $(B, D)$. If we roll to switch, however, it does not matter whether that neighbor is in or out of the independent set: either way it is out at the end of the step and $v$ will be in the independent set. In terms of the bounding state, this allows us to remove the neighbor from $D$ and add $v$ to $B$.

Finally, if there is more than one neighbor of $v$ in $D$, there are cases where we add or do not add regardless of the value of $p_{\text{swap}}$, and so $v$ must be added to $D$ at the next step.

Case III will deal with the subcases where exactly one neighbor of $v$ is in $B$ and we roll to swap. Note that if two or more neighbors of $v$ are in $B$ or if one neighbor of $v$ is in $B$ and we do not roll to swap, the state always remains unchanged and so we leave the bounding state unchanged as well.

If no neighbors of $v$ are in $D$, then the new state of $v$ is completely determined by the neighbors of $v$ that are all known. If at least one neighbor of $v$ is in $D$, then if we do not swap $v$ is always out, but if we do swap then $v$ could be either in or out of the final state, and the neighbor of $v$ could be in or out of the new independent set as well. This means that both must join $D$.

PROOF OF THEOREM 1. The bounding chain detects complete coupling when $|D_t| = 0$ and so we are interested in the behavior of $|D_t|$. In fact, we will upper bound $E[|D_{t+1}||B_t, D_t]$. To make the analysis easier, suppose that we have a line between lines 1 and 2 where we remove node $v$ from $D_t$. We did not do this in our description of the bounding chain simply because it is redundant if $v$ is not in $D_t$ and unnecessary if we end up adding $v$ back to $D_t$. However, we can insert this line without changing the evolution of the bounding chain at all (which case we are in depends only on the neighbors of $v$—not $v$ itself).

This move decreases $|D_t|$ by 1 for all $v \in D_t$, which happens with probability $|D_t|/n$. The remaining lines have a chance of adding one or more nodes back to $D_t$.

If we are in Cases I, IIa or IIIa, we see that nodes are never added to $D_t$. Case IIb is a bad case; $v$ is added back to $D_t$ with probability $(1 - p_{\text{swap}})\lambda/(\lambda + 1)$. Case IIc is a good case; a node other than $v$ is removed from $D_t$ with probability $p_{\text{swap}}\lambda/(\lambda + 1)$. Case IId is another bad case; $v$ is added to $D_t$ with probability $\lambda/(\lambda + 1)$. Case IIIb is the worst case: $|D_t|$ grows by 2 with probability $p_{\text{swap}}\lambda/(\lambda + 1)$. Letting $\alpha = \lambda/(\lambda + 1)$, we have

$$E[|D_{t+1}||B_t, D_t] = |D_t| - \frac{|D_t|}{n} + \frac{N_{01}}{n}(1 - p_{\text{swap}})\alpha - \frac{N_{01}}{n}p_{\text{swap}}\alpha$$
$$+ \frac{N_{02+}}{n}\alpha + 2\frac{N_{11+}}{n}p_{\text{swap}}\alpha,$$

where $N_{01}$ is the number of nodes with $|N_v \cap B| = 0$ and $|N_v \cap D| = 1$, $N_{11+}$ is the number of nodes with $|N_v \cap B| = 1$ and $|N_v \cap D| \geq 1$, and $N_{02+}$ is the number of nodes with $|N_v \cap D| = 0$ and $|N_v \cap B| \geq 2$.

At this point, we take advantage of our ability to choose $p_{\text{swap}}$. Setting $p_{\text{swap}} = 1/4$ means that $1 - 2p_{\text{swap}}$ and $2p_{\text{swap}}$ both are $1/2$. Hence

$$(5) \quad E[|D_{t+1}||B_t, D_t] = |D_t| + \frac{1}{n}\left[-|D_t| + \frac{1}{2}N_{01}\alpha + \frac{1}{2}N_{11+}\alpha + N_{02+}\alpha\right].$$

For a node to count towards $N_{01}$ or $N_{11+}$, it must be adjacent to a node in $D_t$. Similarly, nodes counting towards $N_{02+}$ must adjoin at least two nodes in $D_t$. Since $\Delta$ is the maximum degree of the graph, the number of such nodes is at most $\Delta|D_t|$. Substituting

$$N_{01} + N_{11+} + 2N_{02+} \leq \Delta|D_t|$$

into (5) yields

$$E[|D_{t+1}||B_t, D_t] \leq |D_t|[1 - (1 - \alpha\Delta/2)]$$

and taking expectations of both sides gives

$$E[|D_{t+1}|] \leq E[|D_t|][1 - (1 - \alpha\Delta/2)].$$

A simple induction using $|D_0| = n$ gives us

$$E[|D_t|] \leq n[1 - (1 - \alpha\Delta/2)]^t.$$

As long as $1 - \alpha\Delta/2 > 0$ or, equivalently, $\lambda \leq 2/(\Delta - 2)$, this inequality will be a shrinking upper bound on $E[|D_t|]$. Setting $\beta = [1 - (1 - \alpha\Delta/2)]$ gives us that after $[\log_\beta n] + \theta$ steps, we will have $E[|D_t|] \leq \beta^\theta$ and by Markov's inequality $P(|D_t| \geq 1) \leq \beta^\theta$. □

When $\beta = 1$, we could (using Martingale theory) show that the time until complete coupling is $O(n^2)$. Because such an analysis is of limited utility we omit it here.

In [9], Dyer and Greenhill showed that a simple pairwise coupling will couple in $O(n \ln n)$ time when $p_{\text{swap}} = 1/4$ using the technique of path coupling. This bounding chain not only provides a perfect sampling algorithm, it also serves as an independent proof of their result.

**5. Proper colorings of a graph.** A *proper coloring* of a graph $(V, E)$, where $|V| = n$ is an element $x$ of $\{1, \ldots, k\}^N$ such that for all edges $\{v, w\} \in E$, $x(v) \neq x(w)$.

Consider the Gibbs sampler chain for this problem. Given a proper coloring $x$, let $b_x(v)$ denote the number of different colors not used by neighbors of the node $v$. Then for configurations $x$ and $y$ that differ in color at exactly one node

$$\mathbf{K}(x, y) = \begin{cases} \dfrac{1}{nb_x(v)}, & x, y \text{ differ at node } v, \\ \displaystyle\sum_{v \in N} \dfrac{1}{nb_x(v)}, & x = y, \\ 0, & \text{otherwise.} \end{cases}$$

A bounding chain for this Gibbs sampler was presented independently by the author [16] and Häggström and Nelander [14]. The complete coupling that produces this bounding chain is as follows (here $x$ is the current configuration):

> *The Gibbs sampler for proper colorings*
> 1. **Choose** $v \in_U \{1, \ldots, n\}$, let $N_v$ be the neighbors of $v$
> 2. **Repeat**
> 3.   **Choose** $c \in_U \{1, \ldots, k\}$
> 4. **Until** $c \notin x(N_v)$
> 5. **Let** $x(v) \leftarrow c$

In other words, keep choosing colors for the chosen node uniformly at random from the set of colors until a color is found which does not match any of the neighbors of $v$. This will be the new color for $v$ and is uniform over the set of colors not used by neighbors of $v$.

We use Definition 1 for the bounding chain. Again we choose colors, but this time we might not know whether or not to reject the new color. If $\Delta$ is the maximum degree of the graph, then after choosing $\Delta+1$ different colors we are guaranteed to have picked one that must be accepted. This logic can be codified in the following bounding chain. Let $y$ be the current state.

> *Bounding chain for the Gibbs sampler for proper colorings*
> 1. **Choose** $v \in_U \{1,\ldots,n\}$, let $N_v$ be the neighbors of $v$, let $y(v) \leftarrow \varnothing$
> 2. **Repeat**
> 3.    **Choose** $c \in_U \{1,\ldots,k\}$
> 4.    **If** no neighbor $w$ of $v$ has $y(w) = \{c\}$
> 5.      **Let** $y(v) \leftarrow y(v) \cup \{c\}$
> 6.    **Until** $c \notin \bigcup_{w \in N_v} y(w)$ or $|y(v)| > \Delta$
> 7. **Let** $x(v) \leftarrow c$

THEOREM 2. *When the number of colors $k$ at least $\Delta(\Delta+2)$, let*

$$\beta = 1 - \frac{1 - (\Delta+1)\Delta/[k-\Delta+1]}{n}. \tag{6}$$

*The probability that the bounding chain detects complete coupling in $\log_\beta n + \theta$ steps is at least $1 - \beta^\theta$.*

PROOF. Let $W_t$ denote the number of nodes $v$ with $|Y_t(v)| > 1$. We begin with $W_0 = n$ and wish to find the expected time until $W_t = 0$. Given $Y_t$, we wish to find the expected value of $W_{t+1}$. As in the previous section, if we select a node $v$ with $|Y_t(v)| > 1$, we reduce $W_t$ by 1 at that point.

Now, if we select a color for $v$ that lies within $Y_t(w)$ for some neighbor $w$ of $v$ with $|Y_t(w)| > 1$, that could increase $W_t$ by 1. The probability of this happening is at most $[\sum_{w \in N_v} |Y_t(w)| \mathbb{1}_{|Y_t(v)|>1}]/(k-(\Delta-1))$. (We subtract $\Delta - 1$ from $k$ in the denominator to account for neighbors with known colors that could be preventing some colors from being chosen.) In our bounding chain we always have $|Y_t(w)| \leq \Delta + 1$. Also $\sum_v \sum_{w \in N_v} \mathbb{1}_{|Y_t(v)|>1}$ is just the number of nodes adjacent to nodes with $|Y_t(v)| > 1$ and so is bounded above by $\Delta W_t$. Putting this all together gives us

$$\begin{aligned}
E[W_{t+1}|Y_t] &\leq W_t - \frac{W_t}{n} + \sum_v \frac{1}{n} \sum_{w \in N_v} \frac{|Y_t(w)| \mathbb{1}_{|Y_t(w)|>1}}{k-\Delta+1} \\
&\leq W_t - \frac{W_t}{n} + (\Delta+1)\Delta \frac{W_t}{n} \\
&= W_t \left[1 - \frac{1-(\Delta+1)\Delta/[k-\Delta+1]}{n}\right],
\end{aligned}$$

so

$$E[W_{t+1}] = \beta E[W_t]. \tag{7}$$

Using $W_0 = n$, an induction gives us $E[W_t] = \beta^t n$. As in the previous section, application of Markov's inequality gives us $P(W_t > 1) \leq E[W_t] \leq \beta^t n$ and the result directly follows. □

The best results known for approximately uniform generation of proper colorings is a Markov chain of Vigoda [26] that mixes in polynomial time when $k > 11\Delta/6$. This requires far fewer colors in theory, but when $k$ lies outside this bound the user is given no information. Even when $k < \Delta(\Delta+2)$, a user can still attempt to use the bounding chain approach to generate variates, as (unlike the classical Markov chain method) it is a true algorithm. There is simply no a priori guarantee that the algorithm will terminate in polynomial time.

**6. The antiferromagnetic Potts model.** The Potts model is an extension of the Ising model from statistical physics. Each node of a graph is assigned a color from $\{1, \ldots, k\}$ in a configuration. The energy of configuration $x$ is $H(x) = -\sum_v \sum_{w \in N_v} \mathbb{1}_{x(v)=x(w)}$, where $N_v$ is the set of nodes adjacent to $v$. The probability of choosing a particular configuration is

$$\pi(x) = \frac{1}{Z_T} \exp\left\{\frac{-JH(x)}{T}\right\}. \tag{8}$$

Here $T$ is a positive parameter of the chain known as the temperature, $J$ is either 1 or $-1$, and $Z_T$ is (as with the hard core gas model) known as the partition function. When $J = 1$ more weight is given to configurations with endpoints of edges colored the same way; this is known as the ferromagnetic Potts model. When $J = -1$, configurations where endpoints of edges receive different colors are assigned higher weight: this is the antiferromagnetic Potts model.

In the ferromagnetic case, monotonic chains for alternate formulations of this model exist [12, 23], but for the antiferromagnetic case no chain is known to be monotonic unless the graph is bipartite and $k = 2$. We describe two bounding chains for the antiferromagnetic Potts model based upon the Gibbs sampler. The Gibbs sampler chain chooses a node uniformly at random and then chooses a new color for the node according to $\pi$ conditioned on the value of the rest of the sample.

For the antiferromagnetic Potts model, as the temperature goes to 0 the weight moves towards proper colorings of a graph. For this reason, the uniform distribution over proper colorings of a graph is also known as the antiferromagnetic Potts model at zero temperature.

A bounding chain for the Gibbs sampler can be constructed by noting that each color has a minimum probability of being chosen, regardless of the values of the rest of the sample. This allows the bounding chain to be moved to singleton sets with positive probability. This was the approach used in [16].

Unfortunately, this method has several disadvantages. First, it is time consuming to compute the minimum probabilities for each of the colors. Second, this approach does not result in a provably polynomial time algorithm even when $k$ is large. Intuitively, positive temperature should make it easier to generate random variates, and so we should always be able to generate samples in polynomial time when $k \geq \Delta(\Delta+1)$ for any temperature and with fewer colors as the temperature rises.

We now present a complete coupling for the Gibbs sampler for the antiferromagnetic Potts model. For notational convenience, let $\gamma = \exp(2/T)$.

*The Gibbs sampler for antiferromagnetic Potts model*
1. **Choose** $v \in_U \{1,\ldots,n\}$, let $N_v$ be the neighbors of $v$
2. **Repeat**
3.   **Choose** $c \in_U \{1,\ldots,k\}$
4.   **Choose** $U \in_U [0,1]$
5.   **Let** $a_c$ be the number of neighbors of $v$ with color $c$ in $x$
6. **Until** $U \leq \gamma^{-a_c}$
7. **Let** $x(v) \leftarrow c$

We extend this complete coupling to a bounding chain in the same way as the proper colorings chain: bound $a_c$ between $b_c$ and $d_c$, which are the minimum and maximum number of neighbors of $v$ with color $c$ over the set of $x$ bounded by $y$.

*Bounding chain for the Gibbs sampler for antiferromagnetic Potts model*
1. **Choose** $v \in_U \{1,\ldots,n\}$, let $N_v$ be the neighbors of $v$
2. **Let** $y(v) \leftarrow \varnothing$
3. **Repeat**
4.   **Choose** $c \in_U \{1,\ldots,k\}$
5.   **Choose** $U \in_U [0,1]$
6.   **Let** $b_c$ be the number of $w$ neighboring $v$ with $y(w) = \{c\}$
7.   **Let** $d_c$ be the number of $w$ neighboring $v$ with $c \in y(w)$
8.   **If** $U \leq \gamma^{-b_c}$
9.     **Let** $y(v) \leftarrow y(v) \cup \{c\}$
10. **Until** $U \leq \gamma^{-d_c}$ or $|y(v)| > \Delta$

THEOREM 3. *Let $\tau$ be the first time that the bounding chain detects complete coupling. Then for $k \geq \Delta(\Delta+2)$, and any temperature $T$, or for the case where $\Delta < k < \Delta(\Delta+2)$ and $T$ satisfies $\exp(2/T) < \Delta^2/[\Delta^2+2\Delta-1]$, and finally for the case where $k \leq \Delta$ and $T$ satisfies $\exp(2/T) < \Delta k/[\Delta k-1]$,* *there exists a constant $\beta \in (0,1)$ such that*

$$(9) \qquad P(\tau > (-\ln \beta) n \ln n + \theta) \leq \beta^\theta.$$

PROOF. When $k \geq \Delta(\Delta + 2)$, the same argument used for the proper colorings chain can easily be extended to arbitrary temperatures $T$, whenever we accept when $T = 0$, we will accept when $T > 0$.

Suppose that $\Delta < k < \Delta(\Delta + 2)$. We assume the worst for $Y$, namely, that $|Y_t(v)| = \Delta + 1$ for any $v$ with $|Y_t(v)| > 1$. Let $W_t$ denote the number of nodes with $|Y_t(v)| > 1$. As in the previous section, we will proceed by considering $E[W_{t+1}|Y_t]$. If $|Y_t(w)| > 1$ and $v$ is chosen, $W_t$ is reduced by 1 at the beginning of the step. In the other direction, we will say that node $v$ causes $W_t$ to be increased by 1 if some neighbor $w$ of $v$ is chosen, a color $c$ in $Y_t(v)$ is chosen for that neighbor and we add $c$ to $Y_t(w)$, but we do not exit the repeat loop.

The expected value of $W_{t+1}$ given $Y_t$ will just be $W_t$ plus upper bounds on the expected increase to $W_t$ minus the expected decrease in $W_t$. The expected decrease is just $-W_t/n$. The expected increase caused by particular node $w$ with $|Y_t(w)|$ is always bounded above by $(\Delta + 1)\Delta(1 - \gamma^{-1})/[n(k - \Delta - 1)]$. The factor of $\Delta$ is the maximum number of neighbors for $w$, the $\Delta + 1$ factor the size of $|Y_t(w)|$, the $1 - \gamma^{-1}$ term is the maximum probability that a particular color in $Y_t(w)$ is the culprit for causing $v$ to have a large color set, $1/n$ is the probability of choosing $v$, and $1/(k - \Delta - 1)$ is an upper bound on the probability that the first color for $v$ that we do not reject based on the color of known neighbors is a particular color $c$.

There are $W_t$ such nodes $w$ with $|Y_t(w)| > 1$, so

$$(10) \qquad E[W_{t+1}|Y_t] \leq W_t - \frac{W_t}{n} + W_t\left[\frac{(\Delta + 1)\Delta(1 - \gamma^{-1})}{n(k - \Delta - 1)}\right].$$

Setting $\beta_2 = 1 - (1/n)[(\Delta + 1)\Delta(1 - \gamma^{-1})/(k - \Delta - 1) - 1]$, an induction gives us $E[W_t] \leq \beta_2^t n$ and Markov's inequality, together with the fact that $W_t$ is an integer, allows us to upper bound the probability that $W_t > 0$.

In the final case $k \leq \Delta$ and we can make $|Y_t(w)| = k$. This changes the above argument in two places. First, we replace the factor of $\Delta + 1$ by $k$. Second, we upper bound the probability of choosing a color by 1 rather than $1/(k - \Delta + 1)$. Hence in this case,

$$(11) \qquad E[W_{t+1}|Y_t] \leq W_t - \frac{W_t}{n} + W_t\left[\frac{\Delta k(1 - \gamma^{-1})}{n}\right].$$

Setting $\beta_3 = 1 - (1/n)[\Delta k(1 - \gamma^{-1})]$ gives us $E[W_t] \leq \beta_3^t n$ and the result follows in the same manner as above. □

**7. Sink free orientations of a graph.** Given an undirected graph with node set $N$ and edge set $E$, an *orientation* of the graph is an assignment of direction to each of the individual nodes. A sink free orientation is an orientation such that every node has at least one outgoing edge, and so for this problem $\Omega = \{0,1\}^E$, where 0 and 1 represent each of the two possible directions for each edge. In what follows we will assume without loss of generality the graph is connected, since the number of sink free orientations on an unconnected graph is the product of the number of sink free orientations on each component. Moreover, we will assume that the degree of each node is at least two, since any leaf of the graph fixes the orientation of its adjacent edge outward and can be removed. As with the previous problems, counting the number of sink free orientations of a graph is $\sharp P$-complete [2] and self-reducible so that sample generation leads to an efficient approximation scheme.

The first algorithm for generating almost uniformly from the sink free orientations of a graph (as long as the chain connects the state space) was the basic Gibbs sampler. This was shown to mix in $O([m^3 + mn^3]\ln(1/\varepsilon))$ time by Bubley and Dyer [2] using a simple coupling. Here $\varepsilon$ denotes distance of the variates from uniform distribution in total variation distance.

The first perfect sampling algorithm [16] known for this problem was a bounding chain approach. This method ran in polynomial time, although a technical error in [16] makes the run time theorem proved there invalid. Here we reanalyze that chain with some small improvements to show polynomial running time. Because a new method for perfect sampling known as popping [4] will be faster than any coupling approach of this type, we do not attempt a tight bound on the running time, but rather concentrate on showing that it is polynomial.

The Gibbs sampler chooses a random edge, and then chooses an orientation for the edge uniformly from the set of orientations that do not create a sink. In what follows we will write the orientation of an edge $\{a,b\}$ as $(a,b)$ or $(b,a)$ rather than 0 or 1 to make clear the direction assigned.

---

*Gibbs sampler for sink free orientations of a graph complete coupling*
1.   **Choose** $\{v,w\} \in_U E$
2.   **Choose** $U \in_U [0,1]$
3.   **Case I:** $U < 1/2$ and removing edge $\{w,v\}$ does not make $w$ a sink
4.     **Let** $x(\{v,w\}) \leftarrow (v,w)$
5.   **Case II:** $U > 1/2$ and removing edge $\{w,v\}$ does not make $v$ a sink
6.     **Let** $x(\{v,w\}) \leftarrow (w,v)$

---

We will use Definition 1 for the bounding chain. As previously, we will call an edge $\{v,w\}$ *known* if $|Y(\{v,w\})| = 1$ and *unknown* otherwise. Say that edge $(a,b)$ *leaves* $a$, or its direction from $a$ is *outward*, and it *enters* $b$ or is

directed *inwards* towards $b$. There are several cases to consider in deciding what to do with chosen edge $\{v, w\}$. If there is a known edge leaving $w$ and we choose $(v, w)$, then the edge is known to be oriented in that direction. Similarly, a known edge leaving $v$ allows us to orient the edge $(w, v)$ with certainty.

On the other hand, if all edges adjacent to $w$ are known to enter $w$, then the only orientation possible for the edge is $(w, v)$. Similarly, if all edges adjacent to $v$ are known to enter $v$, the only possible orientation is $(v, w)$. These cases result in the state of the chain (and therefore the bounding chain) remaining unchanged.

Now suppose none of the above cases hold and, in addition, an edge adjacent to $w$ is unknown. Then we do not know if the orientation $(v, w)$ is valid or not: it depends on the state of the unknown edges. Hence if we try to change the orientation to $(v, w)$, we must make $Y(\{v, w\}) = \{(v, w), (w, v)\}$. Similarly, when none of the above cases hold and $v$ is adjacent to an unknown edge, attempted moves to $(w, v)$ result in $\{v, w\}$ becoming unknown in $Y$.

---

*Bounding chain for Gibbs sampler for sink free orientations*
1. **Choose** $\{v, w\} \in_U E$
2. **Choose** $U \in_U [0, 1]$
3. **Case I:** $U < 1/2$ and a known edge is leaving $w$
4.     **Let** $y(\{v, w\}) \leftarrow \{(v, w)\}$
5. **Case II:** $U > 1/2$ and a known edge is leaving $v$
6.     **Let** $y(\{v, w\}) \leftarrow \{(w, v)\}$
7. **Case III:** $U < 1/2$, no known edges leaves $w$, an unknown edge adj. to $w$
8.     **Let** $y(\{v, w\}) \leftarrow \{(v, w), (w, v)\}$
9. **Case IV:** $U > 1/2$, no known edges leaves $v$, an unknown edge adj. to $v$
10.     **Let** $y(\{v, w\}) \leftarrow \{(v, w), (w, v)\}$

---

Unfortunately, the bounding chain built directly from the complete coupling is worthless algorithmically. Suppose all of the edges in the state are unknown. No matter how many steps are taken in the bounding chain, the state of all edges will remain unknown. We need at least one known edge to start the bounding chain.

This is accomplished by using three phases. Phase I consists of a single step in the Markov chain. Given that edge $\{v, w\}$ was chosen and $U < 1/2$, one of two outcomes are possible. Either the edge is oriented $(v, w)$ or the edge is oriented $(w, v)$ and all the other edges adjacent to $w$ are directed into $w$. These are the starting states $Y^1$ and $Y^2$ for two individual bounding chains.

In Phase II we run these bounding chains using independent draws for lines 1 and 2. That is, we make different choices for the edge chosen and $U$ for the $\{Y_t^1\}$ process and the $\{Y_t^2\}$ process. Hopefully we end with each of the two bounding processes (independently) detecting complete coupling.

Finally, in Phase III the two states are run as a regular pairwise coupling just as Bubley and Dyer did in their original paper [2]. If at the end of Phase III, the two states have coupled, then we have achieved complete coupling.

Formally, we have created a complete coupling on the bounded chain with a new state space that is the direct product of the Definition 1 state space and $\{I, II, III\}$ representing which phase we are in. We move from Phase I to Phase II automatically. We have to move from Phase II back to Phase I if either of our bounding process ever reaches the state where all edges are unknown. On the other hand, if both bounding processes detect complete coupling, we get to move on to Phase III.

Since Phase III is the coupling of Bubley and Dyer [3], we use their analysis to conclude that the expected time needed to finish Phase III is $O(m^3 + mn^3)$. Phase I always just takes a single step, and so we are left to analyze Phase II.

Let $W_t$ be the number of unknown edges at time $t$ and consider the expectation of $W_{t+1}$ given $Y_t$. The key thing to note is that an edge can that go from known to unknown must be adjacent to at least 1 unknown edge. Say that edges $\{v, w\}$ and $\{v, z\}$ are *changeable* if one of two conditions hold. In condition (I), $\{v, w\}$ is known to be $(v, w)$ and edge $\{v, z\}$ is unknown, and no other edges are known to leave $v$. By picking $\{v, w\}$ and attempting to direct it towards $v$, we might make this edge unknown and by picking edge $\{v, w\}$ and directing it away from $\{v, w\}$, we might make it known. Hence the contribution to the expected change in $W_t$ coming from this changeable pair is $[1/(2|E|)](-1) + [1/(2|E|)](1) = 0$.

In condition (II), more than one edge is known to be leaving $v$ and an edge $\{v, w\}$ is unknown. If we pick this edge and choose the right direction, this edge becomes known. Hence $W_t$ can only decrease from choosing this type of edge. Since the only edges that change $W_t$ can be partitioned into condition (I) and (II) type edges, the expected value of $W_{t+1}$ given $Y_t$ is at most $W_t$. Thus the stochastic process $W_0, W_1, \ldots$ is a supermartingale.

Suppose that there is a changeable pair. The probability that $W_t$ changes at the next step is at least $1/(2m)$. If no changeable pair exists, then we bide our time, running the chain forward until a pair becomes changeable. Since there are both known and unknown edges in Phase II, there exists a node $a$ adjacent to both known and unknown edges. Let $(b, a)$ be one of the known edges entering $a$ (if any edge leaves $a$, we would have a changeable pair).

Call a node *free* if it is adjacent to at least two known outgoing edges. Note that any edge adjacent to a free node $v$ can be oriented in any direction, moreover, if an edge $(v, w)$ is switched to $(w, v)$, then node $w$ becomes free.

Thus the free nodes take a random walk around the graph. There must be a free node connected to $a$ by known edges. Starting at $a$, take a walk backwards along edges that enter the node we are currently visiting. Either there is an incoming edge to our current node, in which case keep walking, or there are two or more outgoing edges, in which case we have found our free node. The number of nodes in the graph is finite, so eventually we must hit the same node twice, giving us our free node.

On average the free node takes a step at least once every $2m$ moves. After $O(m^3)$ steps it will reach node $b$, at which point $(b,a)$ has a chance of being reoriented $(a,b)$, making the unknown edge adjacent to $a$ together with $(b,a)$ a changeable pair. The value of $W_t$ then has a $1/(2m)$ chance of changing. Hence the expected time between changes to $W_t$ is bounded above by $O(m^5)$.

Now, because $W_t$ is a supermartingale, it has a $1/m$ chance of reaching 0 before it hits $m$, given $W_0 \leq m-1$. We are running the two processes independently and so the probability that each hits 0 before either hits $m$ is at least $1/m^2$. The expected number of steps needed for such a random walk is at most $m^2$ and so altogether, the expected number of steps taken in Phase II is $O(m^2 \cdot m^2 \cdot m^5) = O(m^9)$.

THEOREM 4. *The expected number of steps needed to detect complete coupling for the sink free orientations bounding chain is $O(m^9)$.*

**8. Extensions and open problems.** There are several ways to extend the basic concept of bounding chains. For instance, $V$ or $C$ could be continuous rather than discrete. When $V$ is continuous, the second form of bounding chains is more useful in this regard, but for $C$ continuous the first form is more easily extended [18].

Another extension (both of bounding chains and methodologies such as CFTP) is to deal with more general couplings. All of the complete couplings we discussed in this work are memoryless, dealing with independent draws from $\phi_0, \phi_1, \ldots$. As we saw in the case of the transposition chain, these types of couplings are inherently limited and may not tightly bound the true running time of the chain. To deal with distributions such as uniform over linear extensions of a poset, it is necessary to use more general couplings.

Bounding chains are a straightforward way for determining when a particular complete coupling has in fact completely coupled. Their simplicity is offset somewhat by the fact that the complete coupling time obtained might not tightly bound the mixing time of the chain. Furthermore, no guidance is given for choosing a complete coupling that leads to a good bounding chain. Small changes in a complete coupling can lead to large differences in the behavior of the bounding chain.

However, it is often possible for chains of both theoretical and practical interest to create a bounding chain that can be proven to detect coupling in polynomial time under certain conditions, and which runs quickly under more general situations. While the bounding chain idea is itself simple, it raises a host of interesting theoretical and computational questions about how quickly it moves toward detection of coupling, some of which we were able to answer here, but many of which are still open.

**Acknowledgments.** Much of this work appeared in my Ph.D. thesis, and I would like to thank my advisor David Shmoys for his valuable comments throughout the writing process, as well as Persi Diaconis for many helpful conversations.

I would also like to thank the anonymous referee and Associate Editor for several valuable suggestions.

DEPARTMENT OF MATHEMATICS
AND INSTITUTE OF STATISTICS
AND DECISION SCIENCES
DUKE UNIVERSITY
BOX 90320
DURHAM, NORTH CAROLINA 27713
USA
E-MAIL: mhuber@math.duke.edu
URL: http://www.math.duke.edu/~mhuber